\documentclass[11pt]{amsart}
\usepackage{url}
\begin{document}
\noindent
{\small Topology Atlas Invited Contributions {\bf 6} no.~2 (2001)  5 pp.}
\vspace{\baselineskip}
\title{Partners: Functional Analysis and Topology}
\author{Lawrence Narici}
\address{Department of Mathematics and Computer Science,
St. John's University,
Jamaica, NY 11439, USA}
\email{naricil@stjohns.edu}
\thanks{See {\em Functional Analysis} by G.~Bachman and L.~Narici, 
Dover, Mineola, New York, 2000, a reprint of the the 1966 Academic Press 
book of the same title.
See also the invited contribution {\em What is functional analysis?} 
\cite{Nartaic} by the same author.}
\maketitle

\section*{Introduction}

Functional analysis and topology were born in the first two decades of the
twentieth century and each has greatly influenced the other. Identifying the
dual space---the space of continuous linear functionals---of a normed space
played an especially important role in the formative years of functional
analysis. To further this endeavor, many new kinds (weak, strong, etc.) of
convergence and compactness were introduced . Metric and general topological
spaces evolved in order to provide a framework in which to treat these types
of convergence. As general topology gestated, many concepts were greatly
clarified and simplified. (For example, ``continuous'' meant transforming
convergent sequences into convergent sequences until about 1935.) These
clarifications led to the development of general topological vector spaces
in the 1930's.

\section*{Beginnings}

As set theory developed at the end of the nineteenth century, its paradoxes
revealed that mathematics had a disturbingly shaky foundation. With the aim
of placing set theory in particular and mathematics generally on a firmer
logical pedestal, Hilbert and others looked to Euclidean geometry for a
model \cite{Hil1}. Until that time the  objects of  mathematical attention 
had been quite specific: real numbers, complex numbers, curves,  surfaces. 
Something more general was sought this time. As Hilbert commented:
\begin{quote}
If among my points I consider some systems of things (e.g., the system of
love, law, chimney sweeps \dots) and then accept only my complete axioms 
as the relationships between these things, my theorems (e.g., the 
Pythagorean) are valid for these things also.%
\footnote{He put it another way at a discussion with  some mathematicians 
in the waiting room of the Berlin railway station. He said about his 
geometric axioms ``One must be able at any time to replace {\em 
points, lines and planes} with {\em tables, chairs and beer mugs}.''}
\end{quote}

In other words, ignorance of exactly what the objects were was mandatory.  
``Truth'' was banished, replaced by ``provability''. The new axiomatic
spirit was to consider {\em structures}, arbitrary sets equipped with
operations that obeyed certain rules. This formalist approach dominated
the twentieth century, and is very much still with us.

Various extensions of {\em limit} and {\em continuity} to objects other
than numbers or points has been with us since the 18th century but their
rigorous study---what we might call early ``functional analysis'', in the
sense of analysis on sets whose members were functions---did not begin
until around 1820. Convergence of a sequence of functions meant pointwise
convergence. It was soon realized that imposing more uniformity conditions
was helpful. Stokes and Seidel (1847--8), for example, discovered that
trigonometric series converged with infinitely increasing slowness near a
jump discontinuity and that the discontinuity cannot be enclosed in any
interval in which the convergence is {\em von gleichem Grade} (uniformly
convergent). Heine proved in 1870 that the Fourier series of a piecewise
smooth $2\pi$-periodic function $f$ converges uniformly in any interval
that does not contain a discontinuity of $f$; if $f$ is continuous, then
its Fourier series converges uniformly and absolutely on every closed
interval.  In the presence of uniform convergence, certain attributes
(notably continuity) of each term of a sequence persist to the limit and
series can be integrated term by term. In 1883 Ascoli discovered the
disturbing possibility of a sequence of continuous functions to possess a
discontinuous (pointwise) limit. He found that this behavior disappeared
if the sequence was {\em equicontinuous} \cite{Asc}. These ``uniform''
concepts percolated into analysis generally. In the period 1890--1910,
still other types of convergence of functions were considered such as {\em
relative uniform} convergence and {\em weak} and {\em strong} convergence,
the latter notions being from functional analysis in the modern sense of
the term.

A comprehensive framework for these different kinds of convergence was
evidently desirable. This forced the question: What do you need in order
to talk about convergence? Clearly, a notion of nearness is vital. The
first attempt was Fr\'echet's metric space \cite{Frech1}, then there was
Hausdorff's topological space \cite{Haus}. In the first application of
this set-with-structure approach, Fr\'echet plucked what he deemed to be
the essential properties of distance in the plane (mainly, just the
triangle inequality) and used it to define the metric space. Were the
axioms in use today his only choice as the distillate? Or did he
experiment with weaker requirements? if so, more spaces are brought into
the realm but the number of provable results diminishes. More or stronger
conditions? Then there would be more and better theorems about fewer
things. (Fr\'echet also introduced {\em norm} and the notation $\|\sup 
{\cdot}\|$ for it; the formal definition of normed {\em spaces} was not
given until 1920--1922 by Banach, Hahn and Helly, however.) With the
perspective of the past century, it is well-nigh incredible how much was
deduced from such simple axioms; the same comment of course applies to
topological spaces as well. These two structures alone vindicated faith in
the axiomatic method, albeit with some degenerate cases of
``axiomatics''---defining new things with no other motivation than to
prove theorems about them.

\section*{Geometry and Duality}

In the period 1890--1910, F. Riesz, and E. Schmidt introduced the language
of Euclidean geometry (``orthogonal functions and families'',
``Pythag\-or\-ean theorem'', ``space'', ``dimension'', ``triangle
inequality'', etc.) into Hilbert space. Using Lebesgue's newly minted
integral, Fr\'echet and Riesz commented in 1907 that the space $L_2[a,b]$
of square-integrable functions had a ``geometry'' analogous to that of
``Hilbert space'', i.e., $\ell_2$.

In the same epoch the notions of ``functional'' (a numerical-valued
function whose domain is a set of functions) and ``operator'' (a function
whose domain and range are sets of functions) came into being. This led to
the development of {\em duality} or {\em topological duality}, the study
and use of the {\em continuous dual} $X^\prime$ of all continuous linear
functionals (or ``forms'') on a topological vector space $X$. The
following developments occurred in the period 1900--1918:

\begin{itemize}
\item
(1903; cf.~\cite[pp.~218--227]{BN}) 
In the first formal attempt at describing the topological dual of a normed 
space, Hadamard seeks to characterize the continuous linear functionals on 
the sup-normed space $C[a,b]$ of continuous functions on $[a,b]$. Riesz 
magnificently completes Hadamard's project in 1909; he shows that every 
continuous linear form $f$ on $C[a,b]$ may be written as a Stieltjes 
integral: $f({\cdot}) = \int_{[a,b]} {\cdot}\, dg$, where $g$ is a   
function of bounded variation on $[a,b]$ whose total variation $V(g) = 
\|f\|$.  In today's language we say that $C^\prime[a,b] = NBV[a,b]$, 
where $NBV[a,b]$ denotes the space of normalized functions of bounded 
variation on $[a,b]$ and $=$ signifies surjective norm-isomorphism.
\item
(1907; cf.~\cite[p.~209]{BN}) 
Fr\'echet and Riesz demonstrate that a Hilbert space 
$(X, \langle {\cdot}, {\cdot} \rangle)$ is self-dual: For each continuous 
linear form $f$ on $X$, there is a vector 
$x$ in $X$ such that $f({\cdot}) = \langle {\cdot}, x\rangle$.
\item
(Riesz 1910; cf.~\cite[p.~286]{DS})
The continuous dual $L_p^\prime[a,b]$ of the space $L_p[a,b]$ ($1 < p < 
\infty$) of $p$-th power integrable functions on $[a,b]$ is $L_q[a,b$] 
where $1/p + 1/q = 1$. The analogous result for $\ell_p$ follows in 1913. 
(The discoveries about $L_p[a,b]$ led directly to the general notion of 
normed space.)
\item
(Riesz 1911, inspired by a boundedness notion of Hilbert's; 
cf.~\cite[p.~209]{BN}) 
A linear functional $f$ is continuous if and only if f is ``bounded'' in 
the sense that there is some $M$ such that $|f(x)|$ is at most $M\|x\|$ 
for all $x$.
\item
(Steinhaus 1918; cf.~\cite[p.~289]{DS}) 
$L_1^\prime[a,b] = L_\infty[a,b]$, the space of measurable 
functions $f$ on $[a,b]$ such that, given $f$, there is some $M$ such 
that $|f(x)|$ is at most $M$ for almost all $x$ in $[a,b]$ (f is 
{\em essentially bounded}).
\end{itemize}

To further these investigations in duality, strong use was made of {\em
weak} compactness, {\em weak} and {\em strong} convergence, {\em relative
uniform} convergence, and {\em complete} continuity (mapping weakly
convergent sequences into strongly convergent ones, as Riesz originally
used the term). Some cracks in the metric space approach to provide a
common framework for the various kinds of convergence were visible almost
immediately. Hausdorff's remedy was a more general approach to
``nearness''. Inspired by Hilbert's axioms of open neighborhoods for the
plane, he defined the general topological space in 1914 \cite[Chapters
7--9]{Haus}.

\section{From Metric to Topological Vector Spaces}

With the appearance of Banach's book \cite{Ban1} in 1932, metric
functional analysis (normed, Hilbert and Fr\'echet spaces) had come into
its own. Its stature was elevated when Hilbert space proved to be a
felicitous home for quantum mechanics. But even before 1930 it was known
that pointwise convergence, convergence in measure and compact convergence
eluded description by means of a norm. The treatment of these things in
linear spaces had to await the introduction of locally convex spaces (von
Neumann and Kolmogorov, 1935). It was time for topology to inspire
functional analysis and it was progress in general topology throughout
1930--1940 that enabled the transition from metric linear spaces to
topological vector spaces. With the locally convex space and von Neumann
and Kolmogorov's notion of bounded set (one which is contained in a
sufficiently large scalar multiple of any neighborhood of $0$), duality
theory was transmogrified in the works of Mackey \cite{Mack,Mack1}, and
Grothendieck \cite{Gro,Gro1}. These changes led to Schwartz's theory of
distributions \cite{Schw}.

For further remarks on the developments during this formative period, see
\cite{Dieu1} and the historical remarks in \cite{BTVS}.

\end{document}